# General formulas of global characteristic coefficients of Collatz function

Raouf Rajab
National School of Engineers of Gabes- Tunisia
raouf.rajab@enig.rnu.tn

**Abstract:** The purpose of this paper is to show three general formulas of three global characteristic coefficients of Collatz function. The Collatz function is defined by the following operation on an arbitrary positive integer if N is odd multiply it by 3 and add 1 then the sum obtained is divided by 2, if N is even divide it by 2. Based on the principle, we define the n-order function denoted by $T^n$ such as the different expressions of that function are results of applying Collatz function n times to a natural numbers which are expressed in well determined forms. Based on these expressions, we can characterize that n-order function by three global characteristic coefficients. In the first, we define these three global coefficients. Secondly, we show that each global characteristic coefficient has a general expression as a function of n.

**Key words:** Collatz function; Relative characteristic coefficients; Global characteristic coefficients; Structural order.

**Résumé :** L'objet de cet article est la démonstration des expressions générales de trois coefficients caractéristiques globaux d'une fonction de Collatz d'ordre n. La fonction de Collatz est définie pour tout entier naturel non N comme suit si N est pair on le divise par 2, si non on le multiplie par 3 et on l'ajoute 1 puis la somme obtenue est divisée par 2. En se basant sur ce processus itératif, on définit la fonction de Collatz d'ordre n tel que ces différentes expressions sont obtenues par application de la fonction de Collatz n fois successives sur un ensemble des entiers naturels qui sont exprimés sous des formes bien déterminées. En se basant sur les différentes expressions obtenues, on peut caractériser cette fonction par trois coefficients caractéristiques globaux. Dans un premier temps, on définit ces trois coefficients globaux ensuite on montre que chaque coefficient caractéristique global possède une expression générale qui ne dépend que de n.

**Mots clés :** Fonction de Collatz; Coefficients caractéristiques relatifs; Coefficient caractéristique global; Ordre structurel.

# 1. Introduction et préliminaire

Les suites de Collatz font l'objet des plusieurs travaux de recherche dans l'objectif final est la démonstration de la conjecture de Collatz. Cette conjecture affirme que toutes les suites de Collatz convergent vers le cycle trivial (4,2,1) après un nombre fini des itérations quelque soit le nombre de départ [1]. La particularité principale de cette suite est que l'expression de l'image d'un entier naturel dépend de sa forme de départ (ou de sa parité). Cette particularité engendre un certain nombre des propriétés exceptionnelles et un comportement d'ensemble caractérisé par une régularité structurelle bien déterminée. Cet article fait partie d'un ensemble des travaux qui consistent à déterminer des règles générales ou bien des propriétés générales qui décrivent les comportements de ce genre particulier des suites.

En se basant sur le principe du processus d'itération de la fonction de Collatz (eq 1.1), on peut définir des fonctions de Collatz de différents ordres structurels. Pour une fonction d'ordre structurel n, l'ensemble de départ est subdivisé en $2^n$ sous ensembles ce qui nous permet d'obtenir $2^n$ expressions pour les images. Cette approche fait l'objet des plusieurs études dans l'objectif est de prouver la conjecture de Collatz ou bien de la vérifié pour le plus grand nombre des entiers naturels [2].

Dans cet article, on va définir trois coefficients caractéristiques globaux pour chaque fonction de Collatz d'ordre n quelconque. Ces coefficients sont des fonctions des coefficients qui interviennent dans la détermination des expressions des différentes fonctions de Collatz. Puis on va montrer que chaque coefficient caractéristique global possède une expression générale bien déterminée qui ne dépendent que de l'ordre structurel de la fonction de Collatz considérée.

Dans ce qui suit, on définit quelques éléments de base et on fixe quelques notations concernant les fonctions et les coefficients qui font l'objet de ce travail.

La fonction de Collatz est définie comme suit pour tous les entiers naturels non nuls N [3]:

1.1
$$T^1(N) = \begin{cases} \dfrac{N}{2} & \text{si } N \equiv 0 \pmod 2 \\ \dfrac{3}{2}N + \dfrac{1}{2} & \text{si } N \equiv 1 \pmod 2 \end{cases}$$

On ajoute 1 en exposant pour indiquer qu'elle s'agit d'une fonction de premier ordre structurel.

**Notation 1.1**

Une suite de Collatz de premier terme $P_0$ et de longueur n+1 est notée comme suit:

1.2
$$Sy^1(P_0, n) = (P_0, P_1, \ldots, P_{n-1}, P_n)$$

Avec :

$$T^1(P_k) = P_{k+1} \text{ avec } k \in \mathbb{N}$$

**Notation 1.2**

On adopte la notation suivante pour la $k^{\text{éme}}$ image de $P_0$ obtenue après k itérations successives par application de la fonction $T^1$ sur N:

$$T_k^1(P_0) = P_k$$

Autrement :

$$\underbrace{T^1\left(T^1(\ldots(T^1(T^1(P_0)))\right)}_{\text{K fois}} = T_k^1(P_0)$$

**Notation 1.3**

La suite de Collatz peut s'écrire alors comme suit :

1.3 $$Sy^1(P, n) = (P, T_1^1(P), T_2^1(P), \ldots, T_n^1(P))$$

On désigne par $M_n(P)$ le nombre des entiers impairs contenus dans la suite de Collatz $Sy^1(P, n)$.

Dans le cas général, on peut classer les entiers naturels en $2^n$ sous ensembles au lieu de deux sous ensembles uniquement. Dans ce cas, les entiers naturels sont exprimés sous la forme générale suivante :

1.4 $$P = 2^n a + b_i$$

Avec $b_i \in \{0, 1, 2, 3, \ldots, 2^n - 1\}$ et l'indice $i \in \{1, 2, 3, \ldots, 2^n\}$

$b_i$ est appelé déterminant structurel. L'ensemble des déterminants structurels de la forme d'ordre n est noté $DS(n)$, il correspond à l'ensemble suivant :

1.5 $$DS(n) = \{0, 1, 2, 3, \ldots, 2^n - 1\}$$

On fait appliquer la fonction $T^1$ sur un entier P qui s'écrit sous la forme $2^n a + b_i$ en effectuant n itérations successives, on obtient une image de P qui prend la forme suivante :

1.6 $$T_n^1(2^n a + b_i) = 3^{m(b_i, n)} a + f(b_i, n)$$

Dans ce cas, on parle d'une fonction de Collatz d'ordre structurel n qu'on la note $T^n$ tel que :

$$T^n(2^n a + b_i) = 3^{m(b_i, n)} a + f(b_i, n)$$

La fonction de Collatz d'ordre n notée $T^n$ est une fonction de $\mathbb{N}^*$ dans $\mathbb{N}^*$ tel que tout entier naturel non nul P qui s'écrits sous la forme $2^n a + b$ admet une image obtenue après n itérations successives par application de $T^1$ à l'entier P ceci se traduit par:

1.7 $$T^n(P) = T_n^1(P)$$

Plusieurs chercheurs ont été intéressés par l'étude de cette fonction et des travaux antérieurs ont montrés que $m(b_i, n)$ est un entier naturel correspond au nombre des entiers impairs dans la suite de Collatz suivante $(P, T_1^1(P), T_2^1(P), \ldots, T_{n-1}^1(P))$ alors que le terme $f(b_i, n)$ correspond à la $n^{\text{ième}}$ image de $b_i$ par la fonction $T^1$ de Collatz [2]. On remplace $m(b_i, n)$ et $f(b_i, n)$ par leurs expressions on peut écrire :

1.8 $$T^n(2^n a + b_i) = \begin{cases} a & \text{si } i = 1 \\ 3^{M_{n-1}(b_i)} a + T_n^1(b_i) & \text{si } 2 \leq i \leq 2^n \end{cases}$$

Avec :

- $M_{n-1}(b_i)$ est le nombre des entiers impairs dans la suite de Collatz $Sy^1(b_i, n-1)$.
- $T_n^1(b_i)$ est la $n^{\text{ième}}$ image de $b_i$ par la fonction d'ordre 1 de Collatz.

**Exemple 1.1**

On prend le cas d'un entier qui s'écrit sous la forme d'ordre 4 suivante :

$$N = 2^4 a + 7$$

Donc on a dans le cas n=4 et b=7. On détermine l'expression de son image par deux méthodes différentes. Pour la première, on effectue 4 itérations successives comme suit:

$$T_1^1(N) = 3x2^3 a + 11, T_2^1(N) = 9x2^2 a + 17, T_3^1(N) = 27x2a + 26, T_4^1(N) = 27a + 13$$

L'expression de $T^4(N)$ en fonction de N s'écrit comme ci dessous:

$$T^4(N) = T_4^1(N) = \frac{27}{16}(N - 7) + 13 = \frac{27}{16}N + \frac{19}{16}$$

La deuxième méthode consiste en premier temps à déterminer le nombre des entiers impairs dans la suite suivante :

$$Sy^1(7,3) = (7,11,17,26)$$

Cette suite contient trois entiers impairs, on déduit que :

$$M_3(7) = 3$$

Puis on détermine la valeur de $T_4^1(7)$ :

$$T_4^1(7) = 13$$

On obtient finalement l'expression de $T^4(N)$ comme suit :

$$T^4(N) = \frac{3^3}{2^4}N + \left(13 - \frac{3^3}{2^4}x7\right) = \frac{27}{16}N + \frac{19}{16}$$

Les expressions de la fonction $T^n$ de Collatz

On fait remplacer a par son expression en fonction de P dans l'équation (1.8), on obtient l'expression de $T^n(P)$ en fonction de P comme suit :

1.9 $\quad T^n(P)$

$$= \begin{cases} \dfrac{P}{2^n} & \text{si } P \equiv 0 \pmod{2^n} \\ \dfrac{3^{M_{n-1}(1)}}{2^n}P + \left(T_n^1(1) - \dfrac{3^{M_{n-1}(1)}}{2^n}x1\right) & \text{si } P \equiv 1 \pmod{2^n} \\ \dfrac{3^{M_{n-1}(2)}}{2^n}P + \left(T_n^1(2) - \dfrac{3^{M_{n-1}(2)}}{2^n}x2\right) & \text{si } P \equiv 2 \pmod{2^n} \\ \quad\quad\quad\quad\quad\quad\vdots \\ \dfrac{3^{M_{n-1}(2^n-1)}}{2^n}P + \left(T_n^1(2^n-1) - \dfrac{3^{M_{n-1}(2^n-1)}}{2^n}x(2^n-1)\right) & \text{si } P \equiv 2^n - 1 \pmod{2^n} \end{cases}$$

**Exemple 1.2**

Les différentes expressions de la fonction de Collatz d'ordre 2 notée $T^2$ sont comme suit :

$$T^2(P) = \begin{cases} \dfrac{1}{4}P & \text{si } P \equiv 0 \ (\text{mod} 4) \\ \dfrac{3}{4}P + \dfrac{1}{4} & \text{si } P \equiv 1 \ (\text{mod} 4) \\ \dfrac{3}{4}P + \dfrac{2}{4} & \text{si } P \equiv 2 \ (\text{mod} 4) \\ \dfrac{9}{4}P + \dfrac{5}{4} & \text{si } P \equiv 3 (\text{mod } 4) \end{cases}$$

Par exemple les trois premiers termes de la suite de Collatz de premier terme 7 sont comme suit:

$$T_1^2(7) = 17, T_2^2(7) = T_1^2(17) = 13, T_3^2(7) = 10$$

**Notation    1.4**

On adopte la notation suivante pour une suite de Collatz d'ordre n de premier terme P et de longueur m+1.

$$Sy^n(P, m) = (P, T_1^n(P), T_2^n(P), \ldots, T_{m-1}^n(P), T_m^n(P))$$

Ces suites font l'objet des plusieurs études, elles sont appelées le plus souvent les suites compressées ou aussi les suites extraites.

**Exemple    1.3**

Dans le tableau suivant, on représente quelques termes des trois suites de Collatz de premier terme 7 et de différents ordres structurels.

Tableau 1: Exemples des suites de Collatz de différents ordres structurels

| $T_k^1(7)$ | 11 | 17 | 26 | 13 | 20 | 10 | 5 | 8 | 4 |
|---|---|---|---|---|---|---|---|---|---|
| $T_k^2(7)$ | 11 | 17 | 26 | 13 | 20 | 10 | 5 | 8 | 4 |
| $T_k^3(7)$ | 11 | 17 | 26 | 13 | 20 | 10 | 5 | 8 | 4 |

Ce tableau montre la distribution relative de différents termes des suites de différents ordres structurels. Les suites représentées dans le tableau ci-dessus sont les suivantes:

$$Sy^1(7,9) = (7,11,17,26,13,20,10,5,8,4)$$
$$Sy^2(7,4) = (7,17,13,10,8)$$
$$Sy^3(7,3) = (7,26,10,4)$$

**Notation    1.5**

Pour chaque expression parmi les $2^n$ expressions définissant la fonction $T^n$, on peut distinguer deux coefficients caractéristiques notés comme ci-dessous:

1.10 $$\begin{cases} A_{n,i} = \dfrac{3^{M_{n-1}(b_i)}}{2^n} \\ B_{n,i} = T_n^1(b_i) - \dfrac{3^{M_{n-1}(b_i)}}{2^n} b_i \end{cases}$$

Ces dernières expressions correspondent à des valeurs de $b_i$ allant de 1 à $2^n - 1$ et pour $b_1 = 0$ on sait que la fonction de Collatz n'est pas définie en 0. On peut tirer les valeurs de ces deux coefficients à partir de la première équation du système 1.9 ce qui nous permet d'écrire :

1.11
$$\begin{cases} A_{n,1} = \dfrac{1}{2^n} \\ B_{n,1} = 0 \end{cases}$$

Les coefficients notées $B_{n,i}$ on les appelle les coefficients d'ajustements de la fonction de Collatz d'ordre structurel n puisqu'ils font ajuster les valeurs de $T^n(2^n a + b_i)$ à des entiers naturels.

Les autres coefficients on les appelle coefficients caractéristiques relatifs de la fonction de Collatz considérée.

Comme la fonction de Collatz d'ordre structurel n est définie par $2^n$ expressions donc on peut distinguer $2^n$ coefficients caractéristiques relatifs pour cette fonction qui sont comme suit:

1.12
$$A_{n,1} = \frac{1}{2^n}, A_{n,2} = \frac{3^{M_{n-1}(b_2)}}{2^n}, \ldots, A_{n,2^n} = \frac{3^{M_{n-1}(b_{2^n})}}{2^n}$$

De même on peut déterminer $2^n$ coefficients d'ajustement pour la fonction de Collatz d'ordre n ces coefficients s'écrits comme suit :

1.13
$$B_{n,1} = 0, B_{n,2} = \left(T_n^1(b_2) - \frac{3^{M_{n-1}(b_2)}}{2^n} b_2\right), \ldots, B_{n,2^n} = \left(T_n^1(b_{2^n}) - \frac{3^{M_{n-1}(b_{2^n})}}{2^n} b_{2^n}\right)$$

**Exemple 1.4**

Comme on a quatre expressions différentes pour cette fonction d'ordre 2 donc on peut distinguer quatre coefficients caractéristiques relatifs. Les résultats obtenus pour la fonction $T^2$ sont regroupés dans le tableau ci-dessous :

Tableau 2 : Différentes expressions de la fonction $T^2$ de Collatz

| N | $T^2(N)$ en fonction de a | $T^2(N)$ en fonction de N | $A_{2,i}$ | $B_{2,i}$ |
|---|---|---|---|---|
| 4a | a | $\dfrac{1}{4}N$ | $\dfrac{1}{4}$ | 0 |
| 4a + 1 | 3a + 1 | $\dfrac{3}{4}N + \dfrac{1}{4}$ | $\dfrac{3}{4}$ | $\dfrac{1}{4}$ |
| 4a + 2 | 3a + 2 | $\dfrac{3}{4}N + \dfrac{2}{4}$ | $\dfrac{3}{4}$ | $\dfrac{2}{4}$ |
| 4a + 3 | 9a + 8 | $\dfrac{9}{4}N + \dfrac{5}{4}$ | $\dfrac{9}{4}$ | $\dfrac{5}{4}$ |

L'avant-dernière colonne du tableau ci-dessous contient ces quatre coefficients $A_{2,i}$ et la dernière colonne contient tous les coefficients d'ajustement $B_{2,i}$ de différentes expressions de $T^2$.

**Définition 1.1**

On définit un premier coefficient caractéristiques global de la fonction de Collatz d'ordre n qu'on le note $K_n$ comme le produit de tous les coefficients caractéristiques relatifs. Il s'écrit comme suit :

1.14
$$K_n = \prod_{i=1}^{2^n} A_{n,i}$$

On montre dans cet article que ce coefficient ne dépend que de l'ordre structurel n de la fonction de Collatz considérée et il a pour expression :

$$K_n = \left(\frac{3}{4}\right)^{2^{n-1}n}$$

Dans le cas d'une fonction de Collatz de premier ordre et de second ordre structurel, on vérifie que cette propriété est vraie, en effet le coefficient caractéristique global noté $K_1$ s'écrit comme suit:

$$K_1 = \frac{1}{2} \times \frac{3}{2} = \frac{3}{4} = \left(\frac{3}{4}\right)^{2^{1-1} \times 1}$$

Le premier coefficient caractéristique global relatif à la fonction de Collatz d'ordre 2 s'écrit comme suit :

$$K_2 = \prod_{i=1}^{4} A_{n,i}$$
$$= \frac{1}{4} \times \frac{3}{4} \times \frac{3}{4} \times \frac{9}{4}$$
$$= \frac{81}{256}$$
$$= \left(\frac{3}{4}\right)^4$$
$$= \left(\frac{3}{4}\right)^{2^{2-1} \times 2}$$

**Définition 1.2**

On définit le deuxième coefficient caractéristique global relatif à cette fonction d'ordre n qu'on le note $S_n$ comme la somme de tous les coefficients caractéristiques relatifs donc il est donné par la formule suivante :

1.15
$$S_n = \sum_{i=1}^{2^n} A_{n,i}$$

On montre que ce deuxième coefficient s'écrit sous la forme générale suivante :

$$S_n = 2^n$$

On vérifie que cette relation est vraie pour la fonction de Collatz de premier ordre structurel et de second ordre. On détermine ce deuxième coefficient caractéristique global dans les deux cas n=1 et n=2.

Pour la fonction $T^1$ ce coefficient noté $S_1$ est calculé comme suit :

$$S_1 = \frac{1}{2} + \frac{3}{2} = 2 = 2^1$$

Dans le cas de la fonction $T^2$, le deuxième coefficient caractéristique global s'écrit comme suit :

$$S_2 = \sum_{i=1}^{4} A_{2,i} = \frac{16}{4} = 2^2$$

**Définition     1.3**

On définit le coefficient global d'ajustement comme la somme de tous les coefficients d'ajustement. Il est noté $R_n$ et il a pour expression :

1.16
$$R_n = \sum_{i=1}^{2^n} B_{n,i}$$

On montre que ce coefficient global ne dépend que de n et il s'écrit sous la forme générale suivante:
$$R_n = 2^{n-2}n$$

Dans le cas de la fonction de Collatz de premier ordre structurel, ce coefficient s'écrit comme ci-dessous:

$$R_1 = B_{1,1} + B_{1,2} = 0 + \frac{1}{2} = 2^{1-2}x1$$

La valeur du coefficient d'ajustement global $R_2$ est comme suit :

$$R_2 = \sum_{i=1}^{4} B_{2,i} = 0 + \frac{1}{4} + \frac{2}{4} + \frac{5}{4} = 2 = 2^{2-2}x2$$

On peut conclure que les trois expressions qui on cherche à démontrer pour les trois coefficients globaux déjà définis sont vérifiées pour les deux cas n=1 et n=2.

Les démonstrations de ces expressions générales des coefficients caractéristiques globaux d'une fonction de Collatz d'ordre structurel n font l'objet de cet article. Tout abord, on vérifie qu'elles sont vraies aussi pour le cas n=3, puis on montre par récurrence qu'elles sont vraies pour n'importe quelle fonction de Collatz d'ordre structurel un entier naturel non nul n quelconque.

**Remarque     1.1**

Comme $T_n^1(b) = T^n(b)$ donc on peut remplacer $T_n^1(b)$ par $T^n(b)$ ou vise versa.

**2.     Hypothèses simplificatrices**

On sait que la fonction $T^n$ et l'application $M_n$ ne sont pas définies au point $b_1 = 0$ et les deux coefficients $A_{n,1}$ et $B_{n,1}$ sont déduites à partir de l'équation 1.8 ou bien 1.9. Pour cette raison on fait la distinction entre les deux cas $b_i = 0$ et $b_i \neq 0$ pour les différentes expressions de $A_{n,i}$ et $B_{n,i}$ qui sont données comme suit

$$A_{n,i} = \begin{cases} \dfrac{1}{2^n} & \text{si } i = 1 \\ \dfrac{3^{M_{n-1}(b_i)}}{2^n} & \text{si } 2 \leq i \leq 2^n \end{cases}$$

Les expressions des coefficients d'ajustement sont comme suit :

$$B_{n,i} = \begin{cases} 0 & \text{si } i = 1 \\ T_n^1(b_i) - \dfrac{3^{M_{n-1}(b_i)}}{2^n} b_i & \text{si } 2 \leq i \leq 2^n \end{cases}$$

On remarque que si on pose que :

2.1
$$\begin{cases} M_{n-1}(b_1) = 0 \\ T_n^1(b_1) = 0 \end{cases}$$

On obtient dans ce cas les mêmes valeurs de $A_{n,1}$ et de $B_{n,1}$ comme ci-dessous :

$$\frac{3^{M_{n-1}(b_1)}}{2^n} = \frac{1}{2^n} = A_{n,1}$$

$$T_n^1(b_1) - \frac{3^{M_{n-1}(b_1)}}{2^n} b_1 = 0 = B_{n,1}$$

Ce qui correspond aux valeurs prises par les deux coefficients donc on peut adopter sans aucune ambigüité et pour des raisons de simplification les deux équations hypothétiques (eq2.1). Ces hypothèses ne posent pas aucun problème, elles ne modifient pas le principe de calcul, ni le problème de départ et nous permettent d'obtenir une expression unique pour chaque coefficient. Les deux coefficients peuvent s'écrire sous les formes suivantes pour tout $1 \leq i \leq 2^n$

$$\begin{cases} A_{n,i} = \dfrac{3^{M_{n-1}(b_i)}}{2^n} \\ B_{n,i} = T_n^1(b_i) - \dfrac{3^{M_{n-1}(b_i)}}{2^n} b_i \end{cases}$$

Les expressions des coefficients caractéristiques globaux peuvent être écrites comme suit :

$$\begin{cases} K_n = \displaystyle\prod_{i=1}^{2^n} \dfrac{3^{M_{n-1}(b_i)}}{2^n} \\ S_n = \displaystyle\sum_{i=1}^{2^n} \dfrac{3^{M_{n-1}(b_2)}}{2^n} \\ R_n = \displaystyle\sum_{i=1}^{2^n} \left(T_n^1(b_i) - \dfrac{3^{M_{n-1}(b_i)}}{2^n} b_i\right) \end{cases}$$

## 3. Les coefficients caractéristiques globaux de la fonction de Collatz $T^3$

Dans ce cas, les entiers naturels sont classés en 8 groupes, en effet l'expression générale d'un entier est donnée par la formule générale suivante :

3.1 $\quad\quad\quad\quad\quad\quad\quad N = 2^3 a + b$

Avec a un entier naturel non nul quelconque et b un entier naturel appartenant à l'ensemble $\{0,1,2,3,\ldots,6,7\}$

Les différentes expressions des images correspondent aux différentes formes d'ordre 3 obtenues en effectuant trois itérations à partir de la forme de départ $2^3 a + b$. Ces expressions sont représentées dans le tableau suivant :

Tableau 3 : Différentes expressions de la fonction d'ordre 3 de Collatz

| N | $T^3(N)$ en fonction de a | $T^3(N)$ en fonction de N | $A_{3,i}$ | $B_{3,i}$ |
|---|---|---|---|---|
| $2^3a$ | a | $\frac{1}{8}N$ | $\frac{1}{8}$ | 0 |
| $2^3a+1$ | $3^2a+2$ | $\frac{9}{8}N + \frac{7}{8}$ | $\frac{9}{8}$ | $\frac{7}{8}$ |
| $2^3a+2$ | $3a+1$ | $\frac{3}{8}N + \frac{2}{8}$ | $\frac{3}{8}$ | $\frac{2}{8}$ |
| $2^3a+3$ | $3^2a+4$ | $\frac{9}{8}N + \frac{5}{8}$ | $\frac{9}{8}$ | $\frac{5}{8}$ |
| $2^3a+4$ | $3a+2$ | $\frac{3}{8}N + \frac{4}{8}$ | $\frac{3}{8}$ | $\frac{4}{8}$ |
| $2^3a+5$ | $3a+2$ | $\frac{3}{8}N + \frac{1}{8}$ | $\frac{3}{8}$ | $\frac{1}{8}$ |
| $2^3a+6$ | $3^2a+8$ | $\frac{9}{8}N + \frac{10}{8}$ | $\frac{9}{8}$ | $\frac{10}{8}$ |
| $2^3a+7$ | $3^3a+26$ | $\frac{27}{8}N + \frac{19}{8}$ | $\frac{27}{8}$ | $\frac{19}{8}$ |

L'avant-dernière colonne du tableau ci-dessus contient tous les coefficients caractéristiques relatifs de la fonction $T^3$ donc le coefficient caractéristique global de la fonction d'ordre 3 s'écrit comme suit :

$$K_3 = \prod_{i=1}^{8} A_{3,i}$$
$$= \frac{1}{8} \times \frac{9}{8} \times \frac{3}{8} \times \frac{9}{8} \times \frac{3}{8} \times \frac{3}{8} \times \frac{9}{8} \times \frac{27}{8}$$
$$= \frac{531441}{8^8}$$
$$= (\frac{3}{4})^{12}$$
$$= (\frac{3}{4})^{(2^{3-1} \times 3)}$$

Le deuxième coefficient caractéristique global s'écrit comme suit:

$$S_3 = \sum_{i=1}^{4} A_{3,i} = \frac{64}{8} = 2^3$$

La dernière colonne du tableau ci-dessus contient tous les coefficients d'ajustement de cette fonction ce qui nous permet déduire la valeur de $R_3$ comme suit :

$$R_3 = \sum_{i=1}^{8} B_{3,i} = 0 + \frac{7}{8} + \frac{2}{8} + \frac{5}{8} + \frac{4}{8} + \frac{1}{8} + \frac{10}{8} + \frac{19}{8} = 6 = 2^{3-2} \times 3$$

De même dans ce cas, le calcul prouve que les trois expressions qu'on cherche à démontrer sont vérifiées par les trois coefficients caractéristiques globaux de la fonction de Collatz de troisième ordre structurel.

**4.    Les expressions générales des coefficients caractéristiques globaux**

**Théorème    4.1**

Pour tout entier naturel non nul n, le premier coefficient caractéristique global $K_n$ de la fonction de Collatz d'ordre structurel n à pour expression :

4.1 $$K_n = (\frac{3}{4})^{2^{n-1} \times n}$$

**Théorème    4.2**

Pour tout entier naturel non nul n, le deuxième coefficient caractéristique global noté $S_n$ de la fonction de Collatz d'ordre structurel n à pour expression:

4.2 $$S_n = 2^n$$

**Théorème    4.3**

Pour tout entier naturel non nul n, le coefficient d'ajustement global $R_n$ de la fonction de Collatz d'ordre structurel n à pour expression:

4.3 $$R_n = 2^{n-2} n$$

*Démonstration du théorème 4.1*

On montre par récurrence que cette propriété est vraie pour tout entier naturel non nul n. On sait que la propriété est vraie pour m=1, m=2 et m=3, on suppose que la propriété est vraie pour tout entier m allant de 4 jusqu'a n et on montre que cette propriété est vraie pour m= n+1.

Pour montrer ce théorème on doit déterminer l'expression de l'image d'un entier naturel non nul N par la fonction $T^{n+1}(N)$ et on doit exprimer cet image en fonction de N et $T^n(N)$ ceci ce traduit par

$$T^{n+1}(N) = f(N, T^n(N))$$

L'expression générale de la forme fondamentale d'ordre n d'un entier naturel N est la suivante :

$$N = 2^n a + b_k$$

Avec $b_k \in \{0,1,2,\ldots,2^n - 3, 2^n - 1\}$ et $k \in \{1,2,\ldots,2^n\}$

L'expression de l'image d'un entier de forme $2^n a + b_k$ par la fonction de Collatz d'ordre n s'écrit :

$$T^n(2^n a + b_k) = \frac{3^{M_{n-1}(b_k)}}{2^n} N + T_n^1(b_k) - \frac{3^{M_{n-1}(b_k)}}{2^n} b_k$$

$M_{n-1}(b_k)$ est le nombre des entiers impairs dans la suite $Sy^1(b_k, n-1)$

$T_n^1(b_k)$ est la $n^{\text{ième}}$ image de $b_k$ par la fonction d'ordre 1 de Collatz $T^1$

Le coefficient caractéristique global s'écrit comme suit:

$$K_n = \prod_{k=1}^{2^n} A_{n,k}$$

$$= \frac{3^{Z_n}}{2^{n\,2^n}}$$

Avec :

$$Z_n = \sum_{k=1}^{2^n} M_{n-1}(b_k)$$

En utilisant la relation de récurrence, on peut écrire :

$$K_n = \left(\frac{3}{4}\right)^{2^{n-1} x n}$$

$$= \frac{3^{Z_n}}{2^{n\,2^n}}$$

$$= \frac{3^{Z_n}}{4^{2^{n-1} n}}$$

On déduit que :

4.5
$$Z_n = \sum_{k=1}^{2^n} M_{n-1}(b_k)$$

$$= 2^{n-1}\, n$$

On fait classer les déterminants structurels $(b_k)_{1 \leq k \leq 2^n}$ selon les parités de $T_n^1(b_k)$ donc on procède comme suit

On désigne par p le nombre des déterminants structurels $b_k$ qui vérifient $T_n^1(b_k)$ est pair et par q le nombre des déterminants structurels $b_k$ tel que $T_n^1(b_k)$ est impair donc on peut subdiviser l'ensemble des déterminants structurels d'ordre n en deux sous ensembles selon la parité de $T_n^1(b_k)$ comme ci-dessous :

$$DS(n) = \{b_1, b_2, \dots, b_{2^n}\} \rightarrow \{b_{d_1}, b_{d_2}, \dots, b_{d_p}\}, \{b_{c_1}, b_{c_2}, \dots, b_{c_q}\}$$

Les deux sous ensembles obtenus sont notés comme suit :

$$F = \left\{b_{d_i} \text{ tel que } \forall \text{ entier naturel i } 1 \leq i \leq p \ T_n^1(b_{d_i}) \text{ est pair}\right\}$$

$$G = \left\{b_{c_j} \text{ tel que } \forall \text{ entier naturel j } 1 \leq j \leq q \ T_n^1(b_{c_j}) \text{ est impair }\right\}$$

Pour des raisons de simplification on peut remplacer la notation $b_{d_i}$ par la notation $\alpha_i$ pour tous les éléments de l'ensemble F et de même on remplace la notation $b_{c_j}$ par la notation $\beta_j$ pour tous les éléments de l'ensemble G.

Ce qui nous permet d'écrire :

4.6
$$\begin{cases} F = \{\alpha_1, \alpha_2, \alpha_3, \dots, \alpha_p\} \\ G = \{\beta_1, \beta_2, \beta_3, \dots, \beta_q\} \end{cases}$$

Aussi pour des rasions de simplification, on adopte le changement de notation suivant :

4.7
$$\begin{cases} M_{n-1}(\alpha_i) = u_i \text{ pour tout } 1 \leq i \leq p \\ M_{n-1}(\beta_j) = v_j \text{ pour tout } 1 \leq j \leq q \end{cases}$$

L'expression de l'image d'un entier naturel de forme fondamentale $2^n a + \alpha_i$ s'écrit comme suit :

$$T^n(2^n a + b_{d_i}) = T^n(2^n a + \alpha_i)$$
$$= 3^{u_i} a + T_n^1(\alpha_i)$$

L'expression de l'image d'un entier naturel de forme fondamentale $2^n a + \beta_j$ s'écrit comme suit :

$$T^n(2^n a + b_{c_j}) = T^n(2^n a + \beta_j)$$
$$= 3^{v_j} a + T_n^1(\beta_j)$$

L'expression du coefficient caractéristique global s'écrit:

$$K_n = \left(\frac{3^{u_1}}{2^n} \times \frac{3^{u_2}}{2^n} \times \ldots \times \frac{3^{u_p}}{2^n}\right)\left(\frac{3^{v_1}}{2^n} \times \frac{3^{v_2}}{2^n} \times \ldots \times \frac{3^{v_q}}{2^n}\right)$$

$$= \frac{3^{((u_1+u_2+\cdots+u_p)+(v_1+v_2+\cdots+v_q))}}{2^{(n \times 2^n)}}$$

$$= \frac{3^{Z_n}}{2^{(2^n \times n)}}$$

Ce qui nous permet d'écrire :

$$Z_n = \sum_{i=1}^{p} u_i + \sum_{j=1}^{q} v_j$$
$$= 2^{n-1} n$$

Puisque le nombre total des déterminants structurels pour les formes fondamentales d'ordre n égal à $2^n$ donc évidement on a :

$$p + q = 2^n$$

L'expression de la forme fondamentale d'ordre (n+1) d'un entier naturel N est la suivante :

$$N = 2^{n+1} a + f_r$$

Avec $f_r \in \{0,1,2,\ldots,2^{n+1}-1\}$ et l'indice $r \in \{1,2,\ldots,2^{n+1}\}$

On cherche à exprimé les déterminants structurels $f_r$ relatifs aux formes fondamentales d'ordre (n+1) en fonction des déterminants structurels $b_k$ relatifs aux formes fondamentales d'ordre n.

On sait que l'ensemble des déterminants structurels d'ordre n+1 est le suivant :

$$DS(n+1) = \{0,1,2,\ldots,2^{n+1}-1\}$$

Cet ensemble peut être subdivisé en deux sous ensembles comme suit:

$$\{0,1,2,\ldots,2^{n+1}-1\} \to \{0,1,2,\ldots,2^n-1\}, \{2^n, 2^n+1, \ldots, 2^{n+1}-1\}$$

On pose :

$$\begin{cases} E_1 = \{0,1,2,\ldots,2^n-1\} \\ E_2 = \{2^n, 2^n+1, \ldots, 2^{n+1}-1\} \end{cases}$$

Dans ce cas, on peut écrire :

$$E_1 = \{0,1,2,\ldots,2^n-1\}$$
$$= \{f_1, f_2, \ldots, f_{2^n}\}$$
$$= \{b_1, b_2, \ldots, b_{2^n}\}$$

$$= DS(n)$$
$$E_2 = \{2^n, 2^n + 1, \ldots, 2^{n+1} - 1\}$$
$$= \{f_{(1+2^n)}, f_{(2+2^n)}, \ldots, f_{2^{n+1}}\}$$
$$= \{b_1 + 2^n, b_2 + 2^n, \ldots, b_{2^n} + 2^n\}$$

Ce qui nous permet de déduire les deux relations suivantes pour tout entier naturel k tel que $1 \leq k \leq 2^n$.

4.8
$$\begin{cases} f_k = b_k \\ f_{(k+2^n)} = b_k + 2^n \end{cases}$$

On fait classer les éléments de deux ensembles E1 et E2 selon la parité de $T_n^1(f_r)$ en utilisant le classement déjà établi pour les déterminants structurels $b_k$ de la forme d'ordre n.

Comme on a $\{f_1, f_2, \ldots, f_{2^n}\} = \{b_1, b_2, \ldots, b_{2^n}\}$ donc on peut déduire le classement de l'ensemble $G_1$ comme suit:

$$\{f_1, f_2, \ldots, f_{2^n}\} \rightarrow \{(\alpha_1, \alpha_2, \ldots, \alpha_p), (\beta_1, \beta_2, \ldots, \beta_q)\}$$

De même un classement de l'ensemble $E_2$ selon la parité de $T_n^1(f_r)$ est le suivant :

$$\{f_{1+2^n}, f_{2+2^n}, \ldots, f_{2^{n+1}}\} \rightarrow$$
$$\{(\alpha_1 + 2^n, \alpha_2 + 2^n, \ldots, \alpha_p + 2^n), (\beta_1 + 2^n, \beta_2 + 2^n, \ldots, \beta_q + 2^n)\}$$

On cherche à déterminer les expressions des coefficients caractéristiques relatifs de la fonction d'ordre n+1 de Collatz : $T^{n+1}$

On sait que :

$$T^{n+1}(2^{n+1}a + f_k) = T^{n+1}(2^{n+1}a + b_k)$$
$$= T^1(T^n(2^{n+1}a + b_k))$$
$$= T^1(T^n(2^n(2a) + b_k))$$
$$= T^1(T^n(2^nA + b_k)) \text{ avec } A = 2a$$
$$= T^1\left(3^{M_{n-1}(b_k)} A + T_n^1(b_k)\right)$$
$$= T^1\left(3^{M_{n-1}(b_k)} 2a + T_n^1(b_k)\right)$$

L'expression de $T^1\left(3^{M_{n-1}(b_k)} 2a + T_n^1(b_k)\right)$ dépend de la parité de $T_n^1(b_k)$ comme suit :

$$T^{n+1}(2^{n+1}a + b_k) = \begin{cases} 3^{M_{n-1}(b_k)}a + T_{n+1}^1(b_k) & \text{si } T^n(b_k) \text{ est pair} \\ 3^{M_{n-1}(b_k)+1}a + T_{n+1}^1(b_k) & \text{si non} \end{cases}$$

Équivaut à:

$$T^{n+1}(2^{n+1}a + b_k) = \begin{cases} 3^{M_{n-1}(b_k)}a + \frac{1}{2}T_n^1(b_k) & \text{si } T^n(b_k) \text{ est pair} \\ 3^{M_{n-1}(b_k)+1}a + \frac{3}{2}T_n^1(b_k) + \frac{1}{2} & \text{si non} \end{cases}$$

On peut effectuer les changements suivants pour de raisons de simplification :

-si $T^n(b_k)$ est pair on remplace $T^{n+1}(b_k)$ par $T^{n+1}(\alpha_k)$ et $M_{n-1}(b_k)$ par $u_k$.

-si $T^n(b_k)$ est impair on remplace $T^{n+1}(b_k)$ par $T^{n+1}(\beta_k)$ et $M_{n-1}(b_k)$ par $v_k$.

Ceci nous permet d'écrire :

$$T^{n+1}(2^{n+1}a + b_k) = \begin{cases} 3^{u_k}a + \dfrac{1}{2}T_n^{s1}(\alpha_k) & \text{si } T_n^1(b_k) \text{ est pair} \\ 3^{v_k+1}a + \dfrac{3}{2}T_n^{s1}(\beta_k) + \dfrac{1}{2} & \text{si non} \end{cases}$$

De même on peut écrire:

$$\begin{aligned} T^{n+1}\big(2^{n+1}a + f_{(k+2^n)}\big) &= T^{n+1}(2^{n+1}a + 2^n + b_k) \\ &= T^1\big(T^n(2^{n+1}a + 2^n + b_k)\big) \\ &= T^1\big(T^n(2^n(2a+1) + b_k)\big) \\ &= T^1\big(T^n(2^n A + b_k)\big) \text{ avec } A = 2a + 1 \\ &= T^1\big(3^{M_{n-1}(b_k)} A + T_n^1(b_k)\big) \\ &= T^1\big(3^{M_{n-1}(b_k)} 2a + 3^{M_{n-1}(b_k)} + T_n^1(b_k)\big) \end{aligned}$$

Cette dernière expression dépend de la parité de $T_n^1(b_k)$ comme suit

$$T^{n+1}(2^{n+1}a + 2^n + b_k) = \begin{cases} 3^{M_{n-1}(b_k)}a + \dfrac{1}{2}3^{M_{n-1}(b_k)} + \dfrac{1}{2}T_n^1(b_k) & \text{si } T_n^1(b_k) \text{ est impair} \\ 3^{M_{n-1}(b_k)+1}a + \dfrac{3^{M_{n-1}(b_k)+1}}{2} + \dfrac{3}{2}T_n^1(b_k) + \dfrac{1}{2} & \text{si } T_n^1(b_k) \text{ est pair} \end{cases}$$

De même on peut effectuer les changements suivants selon la parité de $T_n^1(b_k)$

$$T^{n+1}(2^{n+1}a + f_{k+2^n}) = \begin{cases} 3^{v_k}a + \dfrac{1}{2}3^{v_k} + \dfrac{1}{2}T_n^1(\beta_k) & \text{si } T_n^1(b_k) \text{ est impair} \\ 3^{u_k+1}a + \dfrac{3^{u_k+1}}{2} + \dfrac{3}{2}T_n^1(\alpha_k) + \dfrac{1}{2} & \text{si } T_n^1(b_k) \text{ est pair} \end{cases}$$

On fait regrouper les résultats trouvés ci-dessus pour chaque condition comme suit:

-Si $T_n^1(b_k)$ est pair :

$$\begin{cases} T^{n+1}(2^{n+1}a + b_k) = 3^{u_k}a + \dfrac{1}{2}T_n^1(\alpha_k) \\ T^{n+1}(2^{n+1}a + b_{k+2^n}) = 3^{u_k+1}a + \dfrac{3^{u_k+1}}{2} + \dfrac{3}{2}T_n^1(\alpha_k) + \dfrac{1}{2} \end{cases}$$

Dans ce cas, on cherche à exprimer $T^{n+1}(N)$ en fonction de N, on pose:

$$\begin{cases} N_1 = 2^{n+1}a_1 + b_k \\ N_2 = 2^{n+1}a_2 + b_{k+2^n} \end{cases}$$

On peut déduire les expressions suivantes en remplaçant $a_1$ et $a_2$ par leurs expressions respectivement en fonction de $N_1$ et de $N_2$:

4.9 $$\begin{cases} T^{n+1}(N_1) = \dfrac{3^{u_k}}{2^{n+1}}N_1 - \dfrac{3^{u_k}}{2^{n+1}}\alpha_k + \dfrac{1}{2}T_n^1(\alpha_k) \\ T^{n+1}(N_2) = \dfrac{3^{u_k+1}}{2^{n+1}}N_2 - \dfrac{3^{u_k+1}}{2^{n+1}}\alpha_k + \dfrac{3}{2}T_n^1(\alpha_k) + \dfrac{1}{2} \end{cases}$$

Comme le nombre des déterminants structurels qui vérifient la condition $T_n^1(b_k)$ est pair est égal à p donc on peut déduire qu'on a :

-p coefficients caractéristiques relatifs qui s'écrivent sous la forme :

$$\dfrac{3^{u_k}}{2^{n+1}} = \dfrac{1}{2}\dfrac{3^{u_k}}{2^n}$$

- p coefficients caractéristiques relatifs qui s'écrits sous la forme :

$$\frac{3^{(u_k+1)}}{2^{n+1}} = \frac{3}{2}\frac{3^{u_k}}{2^n}$$

-Si $T_n^1(b_k)$ est impair :

$$\begin{cases} T^{n+1}(2^{n+1}a + b_k) = 3^{v_k+1}a + \frac{3}{2}T_n^1(\beta_k) + \frac{1}{2} \\ T^{n+1}(2^{n+1}a + f_{k+2^n}) = 3^{v_k}a + \frac{1}{2}3^{v_k} + \frac{1}{2}T_n^1(\beta_k) \end{cases}$$

Dans ce deuxième cas, l'expression de $T^{n+1}(N)$ en fonction de N s'écrit comme suit :

4.10
$$\begin{cases} T^{n+1}(N_1) = \frac{3^{v_k+1}}{2^{n+1}}N_1 - \frac{3^{v_k+1}}{2^{n+1}}\beta_k + \frac{3}{2}T_n^1(\beta_k) + \frac{1}{2} \\ T^{n+1}(N_2) = \frac{3^{v_k}}{2^{n+1}}N_2 - \frac{3^{v_k}}{2^{n+1}}\beta_k + \frac{1}{2}T_n^1(\beta_k) \end{cases}$$

Comme le nombre des déterminants structurels qui vérifient la condition $T_n^1(b_k)$ est impair est égal à q donc on peut déduire qu'on a :

- q coefficients caractéristiques relatifs qui s'écrits sous la forme :

$$\frac{3^{v_k}}{2^{n+1}} = \frac{1}{2}\frac{3^{v_k}}{2^n}$$

- q coefficients caractéristiques relatifs qui peuvent s'écrire sous la forme :

$$\frac{3^{(v_k+1)}}{2^{n+1}} = \frac{3}{2}\frac{3^{v_k}}{2^n}$$

L'expression du coefficient caractéristique global de la fonction de Collatz d'ordre (n+1) s'écrit comme suit :

$$K_{n+1} = \prod_{r=1}^{2^{n+1}} \frac{3^{M_n(f_r)}}{2^{n+1}}$$

$$= \prod_{k=1}^{2^n} \frac{3^{M_n(f_k)}}{2^{n+1}} \prod_{k=1}^{2^n} \frac{3^{M_n(f_{k+2^n})}}{2^{n+1}}$$

$$= \prod_{i=1}^{p} \frac{3^{u_i}}{2^{n+1}} \prod_{i=1}^{p} \frac{3^{(u_i+1)}}{2^{n+1}} \prod_{j=1}^{q} \frac{3^{v_j}}{2^{n+1}} \prod_{j=1}^{q} \frac{3^{(v_j+1)}}{2^{n+1}}$$

$$= \left(\frac{3^{u_1}}{2^{n+1}} \times \ldots \times \frac{3^{u_p}}{2^{n+1}}\right)\left(\frac{3^{(u_1+1)}}{2^{n+1}} \times \ldots \times \frac{3^{(u_p+1)}}{2^{n+1}}\right)\left(\frac{3^{(v_1+1)}}{2^{n+1}} \times \ldots \times \frac{3^{(v_q+1)}}{2^{n+1}}\right)\left(\frac{3^{v_1}}{2^{n+1}} \times \ldots \times \frac{3^{v_q}}{2^{n+1}}\right)$$

$$= \frac{3^{(u_1+\cdots+u_p)} \times 3^{(u_1+1)+(u_2+1)\ldots+(u_p+1)} \times 3^{(v_1+\cdots+v_q)} \times 3^{(v_1+1)+(v_2+1)\ldots+(v_q+1)}}{(2^{n+1})^{2^{n+1}}}$$

On peut déduire l'expression de $Z_{n+1}$ à partir de cette dernière expression comme suit :

$$Z_{n+1} = \sum_{r=1}^{2^{n+1}} M_n(f_r)$$

$$= \sum_{i=1}^{p} u_i + \sum_{j=1}^{q}(v_j + 1) + \sum_{i=1}^{p}(u_i + 1) + \sum_{j=1}^{q} v_i$$

$$= 2\left(\sum_{i=1}^{p} u_i + \sum_{j=1}^{q} v_j\right) + p + q$$

On sait que :

$$\sum_{i=1}^{p} u_i + \sum_{j=1}^{q} v_j = \sum_{k=1}^{2^n} M_{n-1}(b_k)$$
$$= Z_n$$
$$= 2^{n-1} n$$

De plus on a :

$$p + q = 2^n$$

Par conséquent :

$$Z_{n+1} = 2\left(\sum_{i=1}^{p} u_i + \sum_{i=1}^{q} u_i\right) + p + q$$
$$= 2 \times 2^{n-1} n + 2^n$$
$$= 2^n n + 2^n$$
$$= 2^n(n+1)$$

On déduit que :

$$K_{n+1} = \frac{3^{Z_{n+1}}}{2^{(n+1)2^{n+1}}}$$
$$= \frac{3^{2^n(n+1)}}{4^{(n+1)2^n}}$$
$$= (\frac{3}{4})^{2^n(n+1)}$$

Ce qui implique que la propriété est vraie pour m=n+1 donc on peut conclure que la propriété est vraie pour tout entier naturel non nul n.

### *Démonstration du théorème 4.2*

On montre par récurrence que la propriété est vraie pour tout entier naturel non nul n. on sait d'après ce qui précède que la propriété est vraie pour les trois premiers cas c'est à dire pour n=1, n=2 et n=3.

On suppose que la propriété est vraie pour tout entier naturel m allant de 4 jusqu' à n.

En exploitant les résultats trouvés pour les coefficients caractéristiques relatifs d'ordre (n+1) de la fonction de Collatz d'ordre (n+1). Ces coefficients sont subdivisés comme suit :

- p coefficients qui s'écrits sous la forme :

$$\frac{3^{u_i}}{2^{n+1}}$$

- p coefficients qui s'écrits sous la forme :

$$\frac{3^{u_j+1}}{2^{n+1}}$$

- q coefficients qui s'écrits sous la forme :

$$\frac{3^{v_r}}{2^{n+1}}$$

- q coefficients qui s'écrits sous la forme :

$$\frac{3^{v_j+1}}{2^{n+1}}$$

L'expression du deuxième coefficient caractéristique global relatif à la fonction de Collatz d'ordre (n+1) est comme suit :

$$S_{n+1} = \left(\frac{3^{u_1}}{2^{n+1}} + \cdots + \frac{3^{u_p}}{2^{n+1}}\right) + \left(\frac{3^{(u_1+1)}}{2^{n+1}} + \cdots + \frac{3^{(u_p+1)}}{2^{n+1}}\right) + \left(\frac{3^{(v_1+1)}}{2^{n+1}} + \cdots + \frac{3^{(v_q+1)}}{2^{n+1}}\right) + \left(\frac{3^{v_1}}{2^{n+1}} + \cdots + \frac{3^{v_q}}{2^{n+1}}\right)$$

Par définition, l'expression de $S_n$ est comme suit:

$$S_n = \frac{3^{u_1}}{2^n} + \cdots + \frac{3^{u_p}}{2^n} + \frac{3^{v_1}}{2^n} + \cdots + \frac{3^{v_q}}{2^n}$$

Ce qui nous permet d'écrire :

$$\frac{3^{u_1}}{2^{n+1}} + \cdots + \frac{3^{u_p}}{2^{n+1}} + \frac{3^{v_1}}{2^{n+1}} + \cdots + \frac{3^{v_q}}{2^{n+1}} = \frac{1}{2}\left(\frac{3^{u_1}}{2^n} + \cdots + \frac{3^{u_p}}{2^n} + \frac{3^{v_1}}{2^n} + \cdots + \frac{3^{v_q}}{2^n}\right) = \frac{1}{2}S_n$$

De plus on peut déduire ces deux relations suivantes:

$$\begin{cases} \dfrac{3^{(u_1+1)}}{2^{n+1}} + \cdots + \dfrac{3^{(u_p+1)}}{2^{n+1}} = \dfrac{3}{2}\left(\dfrac{3^{u_1}}{2^n} + \cdots + \dfrac{3^{u_p}}{2^n}\right) \\ \dfrac{3^{(v_1+1)}}{2^{n+1}} + \cdots + \dfrac{3^{(v_q+1)}}{2^{n+1}} = \dfrac{3}{2}\left(\dfrac{3^{v_1}}{2^n} + \cdots + \dfrac{3^{v_q}}{2^n}\right) \end{cases}$$

Ce qui nous permet d'écrire:

$$\left(\frac{3^{(u_1+1)}}{2^{n+1}} + \cdots + \frac{3^{(u_p+1)}}{2^{n+1}}\right) + \left(\frac{3^{(v_1+1)}}{2^{n+1}} + \cdots + \frac{3^{(v_q+1)}}{2^{n+1}}\right)$$
$$= \frac{3}{2}\left(\frac{3^{u_1}}{2^n} + \cdots + \frac{3^{u_p}}{2^n} + \frac{3^{v_1}}{2^n} + \cdots + \frac{3^{v_q}}{2^n}\right)$$
$$= \frac{3}{2}S_n$$

D'après ce qui précède on déduit que :

$$S_{n+1} = \frac{1}{2}S_n + \frac{3}{2}S_n = 2S_n$$

Comme on a supposé que :

$$S_n = 2^n$$

Alors on peut déduire que :

$$S_{n+1} = 2^{n+1}$$

La propriété est vraie pour m=n+1 donc on peut conclure que la propriété est vraie pour tout entier naturel non nul n.

### Démonstration du théorème 4.3

On montre par récurrence que la propriété est vraie pour tout entier naturel non nul n. on sait d'après ce qui précède que la propriété est vraie pour les trois premiers cas c'est à dire pour n=1, n=2 et n=3. On suppose que la propriété est vraie pour tout entier m allant de 4 à n et on montre qu'elle vraie pour m=n+1.

On sait que :

$$R_{n+1} = \sum_{i=1}^{2^{n+1}} B_{n+1,i}$$

Comme on a :

$$B_{n+1,i} = T_{n+1}^1(b_i) - \frac{3^{M_n(b_i)}}{2^{n+1}} b_i$$

On peut écrire alors:

$$R_{n+1} = \sum_{i=1}^{2^n} B_{n+1,i} + \sum_{i=2^n+1}^{2^{n+1}} B_{n+1,i}$$

D'après les équations établies pour les expressions de $T^{n+1}(N)$ en fonction de N et de $T^n(N)$ (eq 4.9 et eq 4.10) lors de la première démonstration on peut déduire qu'on a :

- P coefficients d'ajustement qui ont pour forme :

$$\frac{1}{2} T_n^1(\alpha_k) - \frac{3^{u_k}}{2^{n+1}} \alpha_k$$

- P coefficients d'ajustement qui s'écrits sous la forme :

$$\frac{3}{2} T_n^1(\alpha_k) + \frac{1}{2} - \frac{3^{u_k+1}}{2^{n+1}} \alpha_k$$

- q coefficients d'ajustement qui ont pour forme :

$$\frac{3}{2} T_n^1(\beta_k) + \frac{1}{2} - \frac{3^{v_k+1}}{2^{n+1}} \beta_k$$

- q coefficients d'ajustement qui s'écrits sous la forme :

$$\frac{1}{2} T_n^1(\beta_k) - \frac{3^{v_k}}{2^{n+1}} \beta_k$$

Le coefficient global d'ajustement de la fonction de Collatz d'ordre n+1 s'écrit :

$$R_{n+1} = \sum_{i=1}^{p} \left( \frac{1}{2} T_n^1(\alpha_i) - \frac{3^{u_i}}{2^{n+1}} \alpha_i + \frac{3}{2} T_n^1(\alpha_i) + \frac{1}{2} - \frac{3^{u_i+1}}{2^{n+1}} \alpha_i \right) +$$

$$\sum_{j=1}^{q} (\frac{3}{2} T_n^1(\beta_j) + \frac{1}{2} - \frac{3^{v_j+1}}{2^{n+1}} \beta_j + \frac{1}{2} T_n^1(\beta_j) - \frac{3^{v_j}}{2^{n+1}} \beta_j)$$

On peut simplifier les deux sommes ci-dessus comme suit :

(1) $$\sum_{i=1}^{p} \left( \frac{1}{2} T_n^1(\alpha_i) - \frac{1}{2} \times \frac{3^{u_i}}{2^n} \alpha_i + \frac{3}{2} T_n^1(\alpha_i) + \frac{1}{2} - \frac{3}{2} \times \frac{3^{u_i}}{2^n} \alpha_i \right)$$

$$= \sum_{i=1}^{p} \left( 2 T_n^1(\alpha_i) + \frac{1}{2} - 2 \times \frac{3^{u_i}}{2^n} \alpha_i \right)$$

$$= \frac{p}{2} + 2 \sum_{i=1}^{p} \left( -\frac{3^{u_i}}{2^n} \alpha_i + T_n^1(\alpha_i) \right)$$

(2) $$\sum_{j=1}^{q} \left( \frac{3}{2} T_n^1(\beta_j) + \frac{1}{2} - \frac{3^{v_j+1}}{2^{n+1}} \beta_j + \frac{1}{2} T_n^1(\beta_j) - \frac{3^{v_j}}{2^{n+1}} \beta_j \right)$$

$$= \sum_{j=1}^{q}\left(\frac{3}{2}T_n^1(\beta_j) + \frac{1}{2} - \frac{3}{2}\text{x}\frac{3^{v_j}}{2^n}\beta_j\frac{1}{2}T_n^1(\beta_j) - \frac{1}{2}\text{x}\frac{3^{v_j}}{2^n}\beta_j\right)$$

$$= \sum_{j=1}^{q}\left(2T_n^1(\beta_j) + \frac{1}{2} - 2\text{x}\frac{3^{v_j}}{2^n}\beta_j\right)$$

$$= \frac{q}{2} + 2\sum_{j=1}^{q}(T_n^1(\beta_j) - \frac{3^{v_j}}{2^n}\beta_j)$$

L'expression de $R_{n+1}$ devient après simplification de deux sommes précédentes comme suit :

$$R_{n+1} = \frac{p+q}{2} + 2(\sum_{i=1}^{p}\left(T_n^1(\alpha_i) - \frac{3^{u_i}}{2^n}\alpha_i\right) + \sum_{j=1}^{q}(T_n^1(\beta_j) - \frac{3^{v_j}}{2^n}\beta_j+))$$

On sait que :

$$R_n = \sum_{i=1}^{p}\left(-\frac{3^{u_i}}{2^n}\alpha_i + T_n^1(\alpha_i)\right) + \sum_{j=1}^{q}(-\frac{3^{v_j}}{2^n}\beta_j + T_n^1(\beta_j))$$

$$p + q = 2^n$$

Ceci nous permet déduire la relation de récurrence suivante:

$$R_{n+1} = 2^{n-1} + 2R_n$$

Comme on a supposé que:

$$R_n = 2^{n-2}n$$

Donc on peut conclure que:

$$R_{n+1} = 2^{n-1} + 2^{n-1}n$$
$$= 2^{n-1}(n+1)$$

Alors la propriété est vraie pour m= n+1 donc on peut conclure que la propriété est vraie pour tout entier naturel non nul n.

**Corollaire     4.1**

Le coefficient caractéristique global $K_n$ d'une fonction de Collatz d'ordre structurel n tend vers zéro lorsque n tend vers l'infini:

4.11 $$\lim_{n\to\infty} K_n = 0$$

**Corollaire     4.2**

La moyenne arithmétique des coefficients caractéristiques partiels d'une fonction de Collatz d'ordre structurel n est constante et elle ne dépend pas de l'ordre structurel de la fonction considérée. Cette moyenne est donnée par l'expression suivante :

4.12 $$S_{moy}(n) = \frac{1}{2^n}\sum_{i=1}^{2^n} A_{n,i} = \frac{S_n}{2^n} = 1$$

**Corollaire     4.3**

La moyenne arithmétique des coefficients d'ajustement d'une fonction de Collatz d'ordre structurel n est égale à un quart de l'ordre structurel:

4.13 $$R_{moy}(n) = \frac{1}{2^n} \sum_{i=1}^{2^n} B_{n,i} = \frac{R_n}{2^n} = \frac{n}{4}$$

Ces corollaires sont des conséquences directes des trois théorèmes précédentes.

## 5 Discussion des résultats obtenus

Les expressions générales démontrées pour les trois coefficients caractéristiques globaux traduisent un comportement de groupe des suites de Collatz caractérisé par une régularité de nature structurelle bien déterminée. On peut citer deux remarques qui découlent de ce qui précède dont le premier est en relation avec la conjecture de Collatz.

**Remarque    5.1**

La première corollaire indique que lorsque on fait tendre l'ordre structurel n vers l'infini , le coefficient global tend vers 0, comme ce coefficient est le produit des tous les coefficients caractéristiques relatifs , il ya aussi une proportion de ces coefficients relatifs qui tendent aussi vers 0 et on peut vérifier qu' il existe au moins un coefficient qui vérifient cette limite mais la détermination de la proportion de ces coefficients nécessite une démonstration rigoureuse et on ne peut pas baser uniquement sur cette corollaire ou bien sur les différentes l'expressions établies pour déterminer les différentes proportions. Pour la proportion des coefficients relatifs qui vérifient la condition suivante :

$$\lim_{n \to +\infty} A_{n,i} = 0$$

Comme $T^n(P) = A_{n,i}P + B_{n,i}$ donc ceci est équivalent à les limites suivantes lorsque n tend vers l'infini:

$$A_{n,i} \to 0 \Rightarrow T^n(P) \to B_{\infty,i}$$

La suite de Collatz de premier terme P tend vers son coefficient d'ajustement absolu noté $B_{\infty,i}$ donc on peut conclure qu'une suite de premier terme un entier naturel non nul bien déterminé P que vérifie la conjecture de Collatz doit vérifie les deux conditions suivantes:

$$\begin{cases} \lim_{n \to +\infty} A_{n,i} = 0 \\ \lim_{n \to +\infty} B_{n,i} = 1 \end{cases}$$

**Remarque    5.2**

On peut utiliser les trois coefficients caractéristiques globaux définis dans cet article pour caractériser le comportement d'un ensemble des suites de Collatz qui vérifient un certain nombre de conditions. Dans ce qui suit, on définit cet ensemble des suites de Collatz et on montre qu'elles peuvent être caractérisées par ces différents coefficients étudiés dans cet article.

On considère une suite arithmétique de raison arithmétique 1 et de premier terme 1 et de longueur $2^n$ et qu'on la note comme suit:

$$s(n) = (1, 2, \ldots, 2^n)$$

On fait attribuer à chaque terme la lettre indexée $P_k$ tel que:

$$P_1 = 2^n, P_2 = 1, P_3 = 2, P_4 = 3, \ldots, P_{2^n} = 2^n - 1$$

On remarque que chaque terme de cette suite peut s'écrire sous la forme suivante:

$$P_k = 2^n a_k + b_k$$

Tel que pour tout entier k allant de 2 à $2^n$ on a:

$$\begin{cases} a_k = 0 \\ b_k = k - 1 \end{cases}$$

Et pour k=1, on a:

$$\begin{cases} a_1 = 1 \\ b_1 = 0 \end{cases}$$

On fait correspondre à chaque terme $P_k$ de la suite arithmétique considérée une suite de Collatz de longueur n+1 comme suit:

$$Sy^1(P_k, n) = (P_k, T_1^1(P_k), T_2^1(P_k), \ldots, T_{n-1}^1(P_k), T_n^1(P_k))$$

On obtient un ensemble constitué de $2^n$ suites de Collatz de premier ordre structurel qu'on le note D :

$$D = \{Sy^1(1, n), Sy^1(2, n), \ldots, Sy^1(2^n, n)\}$$

On peut caractériser chaque suite de Collatz de l'ensemble D par deux coefficients, on fait exprimer son dernier terme $T_n^1(P_k)$ en fonction de son premier terme $P_k$:

5.1 $$T_n^1(P_k) = F_n(P_k) P_k + \varphi_n(P_k)$$

-Le premier coefficient noté $F_n(P_k)$ on l'appelle coefficient cumulatif principal.

-Le deuxième terme noté $\varphi_n(P_k)$ on l'appelle coefficient cumulatif secondaire.

On sait que l'expression de l'image d'un entier naturel qui s'écrit sous la forme $2^n a + b_k$ par la fonction de Collatz d'ordre structurel n est comme suit :

$$T^n(P) = A_{n,k} P + B_{n,k}$$

De plus on sait que :

$$T_n^1(P) = T^n(P)$$

On conclu que les deux coefficients caractéristiques de la suite de Collatz de premier terme $P_k$ et de longueur n+1 correspond aux deux coefficients caractéristiques relatifs de la fonction de Collatz d'ordre n, ceci se traduit par:

5.2 $$\begin{cases} A_{n,k} = F_n(P_k) \\ B_{n,k} = \varphi_n(P_k) \end{cases}$$

Avec $1 \leq k \leq 2^n$

Les expressions des trois coefficients caractéristiques globaux de la fonction de Collatz peuvent être écrites en fonction des coefficients cumulatifs des suites de Collatz de l'ensemble D comme ci-dessous:

5.3 $$\begin{cases} K_n = \prod_{k=1}^{2^n} F_n(P_k) \\ S_n = \sum_{k=1}^{2^n} F_n(P_k) \\ R_n = \sum_{k=1}^{2^n} \varphi_n(P_k) \end{cases}$$

La représentation sous forme d'un tableau de l'ensemble D correspond à un tableau de n+1 colonnes et de $2^n$ lignes tel que la ligne du rang k du tableau considéré renferme tous les termes de la suite $Sy^1(P_k, n)$. On peut l'appeler le tableau complet de Collatz, Il est représenté ci dessus:

| | $Sy^1(P_k, n)$ | | | | | Coefficients relatifs | |
|---|---|---|---|---|---|---|---|
| $P_k$ | $T_1^1(P_k)$ | $T_2^1(P_k)$ | | $T_{n-1}^1(P_k)$ | $T_n^1(P_k)$ | $A_{n,k}$ | $B_{n,k}$ |
| 1 | 2 | 1 | | | | | |
| 2 | 1 | 2 | | | | | |
| 3 | 5 | 8 | | | | | |
| 4 | 2 | 1 | | | | | |
| | | | | | | | |
| | | | | | | | |
| $2^n$ | $2^{n-1}$ | $2^{n-2}$ | | | | | |

Fig 1: Tableau complet de Collatz

Par conclusion :

-Le coefficient $K_n$ correspond au produit de tous les coefficients cumulatifs $F_n(P_k)$ de toutes les suites du tableau considéré.

-Le coefficient $S_n$ égal à la somme de tous les coefficients cumulatifs $F_n(P_k)$ de toutes les suites du l'ensemble considéré.

-Le coefficient $R_n$ égal à la somme de tous les coefficients cumulatifs secondaires $\varphi_n(P_k)$ de toutes les suites du tableau ci-dessus.

Des études sur les comportements structurels des suites de Collatz peuvent montrer l'utilité et l'importance des ces coefficients caractéristiques globaux dans la détermination d'autres propriétés intéressantes qui nous permettent de bien comprendre ce genre particulier des suites numériques. A titre d'exemple si on fait classer les coefficients caractéristiques relatifs de la fonction de Collatz, on fait les comparer par rapport à 1, on peut distinguer deux catégories:

-Les coefficients $A^+$ qui vérifient la condition suivante:

$$\frac{3^{M_{n-1}(b_i)}}{2^n} > 1$$

-Les coefficients $A^-$ qui vérifient la condition suivante:

$$\frac{3^{M_{n-1}(b_i)}}{2^n} < 1$$

On peut montrer en exploitant les expressions établies dans cet article que la proportion des coefficients $A^+$ temps vers 0 lorsque n tend vers l'infini.

# RÉFÉRENCES